\documentclass{compositio}

\usepackage{amsmath}
\usepackage{amsthm,amsfonts,amssymb,mathrsfs}
\usepackage{epic,eepic}
\usepackage{yfonts}
\usepackage{paralist,enumerate}
\usepackage[all]{xy}

\newtheorem{theorem}{Theorem}

\newtheorem{corollary}[theorem]{Corollary}

\newtheorem*{Brunn-Minkowski}{Brunn-Minkowski Inequality}
\newtheorem*{proposition*}{Proposition}
\newtheorem*{conjecture*}{Conjecture}
\newtheorem*{question*}{Question}

\newtheorem*{problem*}{Problem}
\newtheorem*{theorem*}{Theorem}

\theoremstyle{definition}

\newtheorem*{definition*}{Definition}
\newtheorem*{Var}{Varchenko's conjecture}

\theoremstyle{remark}
\newtheorem{remark}[theorem]{Remark}

\newtheorem*{example*}{Example}
\newtheorem*{remark*}{Remark}

\usepackage{ifpdf}
\ifpdf
  \usepackage{graphicx}
  \usepackage{xcolor}
\else
  \usepackage[dvipdfmx]{graphicx}
  \usepackage[dvipdfmx]{xcolor}
\fi

\begin{document}

\title[$h$-vectors of matroids and logarithmic concavity]{$h$-vectors of matroids and logarithmic concavity}
\author{June Huh}
\email{junehuh@umich.edu}
\address{Department of Mathematics, University of Michigan\\
Ann Arbor, MI 48109 \\ USA}
\begin{abstract}
Let $M$ be a matroid on $E$, representable over a field of characteristic zero. We show that $h$-vectors of the following simplicial complexes are log-concave:
\begin{enumerate}[1.]
\item The matroid complex of independent subsets of $E$.
\item The broken circuit complex of $M$ relative to an ordering of $E$.
\end{enumerate}
The first implies a conjecture of Colbourn on the reliability polynomial of a graph, and the second implies a conjecture of Hoggar on the chromatic polynomial of a graph.
The proof is based on the geometric formula for the characteristic polynomial of Denham, Garrousian, and Schulze.
\end{abstract}

\maketitle 

\section{Introduction and results}

A sequence $e_0,e_1,\ldots,e_n$ of integers is said to be 
\emph{log-concave} if for all $0< i< n$,
\[
e_{i-1}e_{i+1}\le e_i^2,
\]
and is said to have \emph{no internal zeros} if there do not exist $i < j <  k$ satisfying 
\[
e_i \neq 0, \quad e_j=0, \quad e_k\neq 0.
\] 
Empirical evidence has suggested that many important enumerative sequences are log-concave, but proving the log-concavity can sometimes be a non-trivial task. See \cite{Brenti,Stanley,Stanley00} for a wealth of examples arising from algebra, geometry, and combinatorics. The purpose of this paper is to demonstrate the use of an algebro-geometric tool to the log-concavity problems.

Let $X$ be a complex algebraic variety. A subvariety of $X$ is an irreducible closed algebraic subset of $X$. 
If $V$ is a subvariety of $X$, then the top dimensional homology group $H_{2\dim(V)}(V;\mathbb{Z}) \simeq \mathbb{Z}$ has a canonical generator, and the closed embedding of $V$ in $X$ determines a homomorphism
\[
H_{2\dim(V)}(V;\mathbb{Z}) \longrightarrow H_{2\dim(V)}(X;\mathbb{Z}).
\]
The image of the generator is called the \emph{fundamental class} of $V$ in $X$, denoted $[V]$. A homology class in $H_*(X;\mathbb{Z})$ is said to be \emph{representable} if it is the fundamental class of a subvariety. 

Hartshorne asks in \cite[Question 1.3]{HartshorneSurvey} 
which even dimensional homology classes of $X$ are representable by a smooth subvariety. Although the question is exceedingly difficult in general, it has a simple partial answer when $X$ is the product of complex projective spaces  $\mathbb{P}^m \times \mathbb{P}^n$. Note in this case that the $2k$-dimensional homology group of $X$ is freely generated by the classes of subvarieties of the form $\mathbb{P}^{k-i} \times \mathbb{P}^{i}$. 


Representable homology classes of $\mathbb{P}^m \times \mathbb{P}^n$ can be characterized numerically as follows \cite[Theorem 20]{Huh}.

\begin{theorem}\label{LC}
Write $\xi \in H_{2k}(\mathbb{P}^m \times \mathbb{P}^n;\mathbb{Z})$ as the integral linear combination
\[
\xi =\sum_{i} e_i \big[\mathbb{P}^{k-i} \times \mathbb{P}^{i}\big].
\]
\begin{enumerate}
\item If $\xi$ is an integer multiple of either
\[
\big[\mathbb{P}^m \times \mathbb{P}^n\big], \big[\mathbb{P}^m \times \mathbb{P}^0\big], \big[\mathbb{P}^0 \times \mathbb{P}^n\big], \big[\mathbb{P}^0 \times \mathbb{P}^0\big],
\]
then $\xi$ is representable if and only if the integer is $1$.
\item If otherwise, some positive integer multiple of $\xi$ is representable if and only if the $e_i$ form a nonzero log-concave sequence of
nonnegative integers with no internal zeros.
\end{enumerate}
\end{theorem}

In short, subvarieties of $\mathbb{P}^m \times \mathbb{P}^n$ correspond to log-concave sequences of nonnegative integers with no internal zeros. Therefore, when trying to prove the log-concavity of a sequence, it is reasonable to look for a subvariety of $\mathbb{P}^m \times \mathbb{P}^n$ which witnesses this property. We demonstrate this method by proving the log-concavity of $h$-vectors of two simplicial complexes associated to a matroid, when the matroid is representable over a field of characteristic zero. Other illustrations can be found in \cite{Huh,Huh-Katz,Lenz}.

In order to fix notations, we recall from \cite{Bjorner3} some basic definitions on simplicial complexes associated to a matroid. We use Oxley's book as our basic reference on matroid theory \cite{Oxley}. 

Let $\Delta$ be an abstract simplicial complex of dimension $r$. The \emph{$f$-vector} of $\Delta$ is a sequence of integers $f_0,f_1,\ldots,f_{r+1}$, where
\[
f_i = \big(\text{the number of $(i-1)$-dimensional faces of $\Delta$}\big).
\]
For example, $f_0$ is one, $f_1$ is the number of vertices of $\Delta$, and $f_{r+1}$ is the number of facets of $\Delta$. The \emph{$h$-vector} of $\Delta$ is defined from the $f$-vector by the polynomial identity
\[
\sum_{i=0}^{r+1} f_i (q-1)^{r+1-i}= \sum_{i=0}^{r+1} h_i q^{r+1-i}.
\]
When there is a need for clarification, we write the coefficients by $f_i(\Delta)$ and $h_i(\Delta)$ respectively.

Let $M$ be a matroid of rank $r+1$ on an ordered set $E$ of cardinality $n+1$. We are interested in the $h$-vectors of the following simplicial complexes associated to $M$:
\begin{enumerate}[1.]
\item The matroid complex $\text{IN}(M)$, the collection of subsets of $E$ which are independent in $M$.
\item The broken circuit complex $\text{BC}(M)$, the collection of subsets of $E$ which do not contain any broken circuit of $M$. 
\end{enumerate}
Recall that a \emph{broken circuit} is a subset of $E$ obtained from a circuit of $M$ by deleting the least element relative to the ordering of $E$.
We note that the isomorphism type of the broken circuit complex does depend on the ordering of $E$. However, the results of this paper will be independent of the ordering of $E$.

\begin{remark}
A pure $r$-dimensional simplicial complex is said to be \emph{shellable} if there is an ordering of its facets
such that each facet intersects the complex generated by its predecessors in a pure $(r-1)$-dimensional complex.
$\text{IN}(M)$ and $\text{BC}(M)$ are pure of dimension $r$, and are shellable. As a consequence, the $h$-vectors of both complexes consist of nonnegative integers \cite{Bjorner3}. This nonnegativity is recovered in Theorem \ref{main} below.
\end{remark}


Dawson conjectured that the $h$-vector of a matroid complex is a log-concave sequence \cite[Conjecture 2.5]{Dawson}. Our main result verifies this conjecture for matroids representable over a field of characteristic zero.

\begin{theorem}\label{main}
Let $M$ be a matroid representable over a field of characteristic zero.
\begin{enumerate}
\item  The $h$-vector of the matroid complex of $M$ is a log-concave sequence of nonnegative integers with no internal zeros.
\item The $h$-vector of the broken circuit complex of $M$ is a log-concave sequence of nonnegative integers with no internal zeros.
\end{enumerate}
\end{theorem}

Indeed, as we explain in the following section, there is a subvariety of a product of projective spaces which witnesses the validity of Theorem \ref{main}.


It can be shown that the log-concavity of the $h$-vector implies the strict log-concavity of the $f$-vector:
\[
f_{i-1}f_{i+1}<f_i^2, \qquad i=1,2,\ldots,r.
\]
See  \cite[Lemma 5.1]{Lenz}. Therefore Theorem \ref{main} implies that the two $f$-vectors associated to $M$ are strictly log-concave.
The first statement of the following corollary recovers \cite[Theorem 1.1]{Lenz}.

\begin{corollary}\label{f-vector}
Let $M$ be a matroid representable over a field of characteristic zero.
\begin{enumerate}
\item  The $f$-vector of the matroid complex of $M$ is a strictly log-concave sequence of nonnegative integers with no internal zeros.
\item The $f$-vector of the broken circuit complex of $M$ is a strictly log-concave sequence of nonnegative integers with no internal zeros.
\end{enumerate}
\end{corollary}

The main special cases of Theorem \ref{main} and Corollary \ref{f-vector} are treated in the following subsections.

\begin{remark}
A pure simplicial complex is a matroid complex if and only if every ordering of the vertices induces a shelling \cite[Theorem 7.3.4]{Bjorner3}. 
In view of this characterization of matroids, one should contrast Theorem \ref{main} with examples of other `nice' shellable simplicial complexes whose $f$-vector and $h$-vector fail to be log-concave. In fact, the unimodality of the $f$-vector already fails for simplicial polytopes in dimension $\ge 20$ \cite{Billera-Lee,Bjorner1}. 

These shellable simplicial complexes led to suspect that various log-concavity conjectures on matroids might not be true in general \cite{Stanley00,Wagner}.
Theorem \ref{main} shows that there is a qualitative difference between the $h$-vectors of
\begin{enumerate}[1.]
\item matroid complexes and other shellable simplicial complexes, and/or 
\item matroids representable over a field and matroids in general.
\end{enumerate}
See \cite{Stan} and \cite{StanCM} for characterizations of $h$-vectors of simplicial polytopes and, respectively, shellable simplicial complexes in general. We note that the method of the present paper to prove the log-concavity crucially depends on the assumption that the matroid is representable over a field.
\end{remark}

\subsection{The reliability polynomial of a graph}

The \emph{reliability} of a connected graph $G$ is the probability that the graph remains connected when each edge is independently removed with the same probability $1-p$. If the graph has $e$ edges and $v$ vertices, then the reliability of $G$ is the polynomial
\[
\text{Rel}_G(p) = \sum_{i=0}^{e-v+1} f_i \hspace{0.5mm} p^{e-i}(1-p)^i,
\]
where $f_i$ is the number of cardinality $i$ sets of edges whose removal does not disconnect $G$. For example, $f_0$ is one, $f_1$ is the number of edges of $G$ that are not isthmuses, and $f_{e-v+1}$ is the number of spanning trees of $G$. The \emph{$h$-sequence} of the reliability polynomial is the sequence $h_i$ defined by the expression
\[
\text{Rel}_G(p) = p^{v-1} \Bigg(\sum_{i=0}^{e-v+1} h_i \hspace{0.5mm} (1-p)^i\Bigg).
\] 
In other words, the $h$-sequence is the $h$-vector of the matroid complex of the cocycle matroid of $G$.
Since the cocycle matroid of a graph is representable over every field, Theorem \ref{main} confirms a conjecture of Colbourn
that the $h$-sequence of the reliability polynomial of a graph is log-concave \cite{Colbourn}.

\begin{corollary}\label{Colbourn's conjecture}
The $h$-sequence of the reliability polynomial of a connected graph is a log-concave sequence of nonnegative integers with no internal zeros.
\end{corollary}

It has been suggested that Corollary \ref{Colbourn's conjecture} has practical applications in combinatorial reliability theory \cite{Brown-Colbourn}.

\subsection{The chromatic polynomial of a graph}

The \emph{chromatic polynomial} of a graph $G$ is the polynomial defined by
\[
\chi_G(q) = (\text{the number of proper colorings of $G$ using $q$ colors}). 
\]
The chromatic polynomial depends only on the cycle matroid of the graph, up to a factor of the form $q^c$. More precisely, the absolute value of the $i$-th coefficient of the chromatic polynomial is the number of cardinality $i$ sets of edges which contain no broken circuit \cite{Whitney}.
Since the cycle matroid of a graph is representable over every field, Corollary \ref{f-vector} confirms a conjecture of Hoggar that the coefficients of the chromatic polynomial of a graph form a strictly log-concave sequence \cite{Hoggar}. 
 
\begin{corollary}\label{sub}
The coefficients of the chromatic polynomial of a graph form a sign-alternating strictly log-concave
sequence of integers with no internal zeros.
\end{corollary}

Corollary \ref{sub} has been previously verified for all graphs with $\le 11$ vertices \cite{Lundow-Markstrom}.



\section{Proof of Theorem \ref{main}}

We shall assume familiarity with the M\"obius function $\mu(x,y)$ of the lattice of flats $\mathscr{L}_M$. For this and more, we refer to \cite{Aigner,Zaslavsky}. An important role will be played by the \emph{characteristic polynomial} $\chi_M(q)$. For a loopless matroid $M$, the characteristic polynomial is defined from $\mathscr{L}_M$ by the formula
\[
\chi_M(q) = \sum_{x \in \mathscr{L}_M} \mu(\varnothing, x) q^{r+1 - \text{rank} (x)}=\sum_{i=0}^{r+1} (-1)^i w_i q^{r+1-i}.
\]
If $M$ has a loop, then $\chi_M(q)$ is defined to be the zero polynomial. The nonnegative integers $w_i$ are called the \emph{Whitney numbers of the first kind}. The characteristic polynomial is always divisible by $q-1$, defining the \emph{reduced characteristic polynomial}
\[
\overline{\chi_M}(q)=\chi_M(q)/(q-1).
\]

\subsection{Brylawski's theorem I}

We need to quote a few results from Brylawski's analysis on the broken circuit complex \cite{Brylawski}. 
The first of these says that the Whitney number $w_i$ is the number of cardinality $i$ subsets of $E$ which contain no broken circuit relative to any fixed ordering of $E$ \cite[Theorem 3.3]{Brylawski}.
This observation goes back to Hassler Whitney, who stated it for graphs \cite{Whitney}.

Fix an ordering of $E$, and let $0$ be the smallest element of $E$. We write $\overline{\text{BC}}(M)$ for the \emph{reduced broken circuit complex} of $M$, the family of all subsets of $E \setminus \{0\}$ that do not contain any broken circuit of $M$. Since the broken circuit complex is the cone over $\overline{\text{BC}}(M)$ with apex $0$, the above quoted fact says that
\[
\overline{\chi_M}(q)= \sum_{i=0}^r (-1)^if_i\big(\overline{\text{BC}}(M)\big) q^{r-i}.
\]
In terms of the $h$-vector, we have
\[
\overline{\chi_M}(q+1)= \sum_{i=0}^r (-1)^{i} h_i\big(\overline{\text{BC}}(M)\big) q^{r-i}=\sum_{i=0}^r (-1)^{i} h_i\big(\text{BC}(M)\big) q^{r-i}, \quad h_{r+1}\big(\text{BC}(M)\big)=0.
\]
Therefore the second assertion of Theorem \ref{main} is equivalent to the statement that the coefficients of $\overline{\chi_M}(q+1)$ form a sign-alternating log-concave sequence with no internal zeros. 

\subsection{Brylawski's theorem II}

We show that the first assertion of Theorem \ref{main} is implied by the second. This follows from the fact that the matroid complex of $M$ is the reduced broken circuit complex of the free dual extension of $M$ \cite[Theorem 4.2]{Brylawski}. We note that not every reduced broken circuit complex can be realized as a matroid complex \cite[Remark 4.3]{Brylawski}. The second assertion of Theorem \ref{main} is strictly stronger than the first in this sense.


Recall that the \emph{free dual extension} of $M$ is defined by taking the dual of $M$, placing a new element $p$ in general position (taking the free extension), and again taking the dual. In symbols,
\[
M \times p := (M^* + p)^*.
\]
If $M$ is representable over a field, then $M \times p$ is representable over some finite extension of the same field. 
Choose an ordering of $E \cup \{p\}$ such that $p$ is smaller than any other element.
Then, with respect to the chosen ordering,
\[
\text{IN}(M)=\overline{\text{BC}}(M \times  p).
\]
For more details on the free dual extension, see \cite{Brylawski,BrylawskiBook,Lenz}.

\subsection{Reduction to simple matroids}

A standard argument shows that it is enough to prove the assertion on $\overline{\chi_M}(q+1)$ when $M$ is simple:

\begin{enumerate}[1.]
\item If $M$ has a loop, then the reduced characteristic polynomial of $M$ is zero, so there is nothing to show in this case.
\item If $M$ is loopless but has parallel elements, replace $M$ by its \emph{simplification} $\overline{M}$ as defined in \cite[Section 1.7]{Oxley}. Then 
the reduced characteristic polynomials of $M$ and $\overline{M}$ coincide 
because $\mathscr{L}_M \simeq \mathscr{L}_{\overline{M}}$. 
\end{enumerate}
Hereafter $M$ is assumed to be simple of rank $r+1$ with $n+1$ elements, representable over a field of characteristic zero.

\subsection{Reduction to complex hyperplane arrangements}

We reduce the main assertion to the case of essential arrangements of affine hyperplanes.
We use the book of Orlik and Terao as our basic reference in hyperplane arrangements \cite{Orlik-Terao}.

Note that the condition of representability for matroids of given rank and given number of elements can be expressed in a first-order sentence in the language of fields.
Since the theory of algebraically closed fields of characteristic zero is complete \cite[Corollary 3.2.3]{Marker},
a matroid representable over a field of characteristic zero is in fact representable over $\mathbb{C}$.

Let $\widetilde{\mathcal{A}}$ be a central arrangement of $n+1$ distinct hyperplanes in $\mathbb{C}^{r+1}$ representing $M$. 
This means that there is a bijective correspondence between $E$ and the set of hyperplanes of $\widetilde{\mathcal{A}}$ which identifies the geometric lattice $\mathscr{L}_M$ with the lattice of flats of $\widetilde{\mathcal{A}}$.
Choose any one hyperplane from the projectivization of $\widetilde{\mathcal{A}}$ in $\mathbb{P}^r$. The \emph{decone} of the central arrangement, denoted $\mathcal{A}$, is the essential arrangement of $n$ hyperplanes in $\mathbb{C}^r$ obtained by declaring the chosen hyperplane to be the hyperplane at infinity.
If $\chi_{\mathcal{A}}(q)$ is the characteristic polynomial of the decone, then 
\[
\chi_{\mathcal{A}}(q)=\overline{\chi_M}(q).
\]
Therefore it suffices to prove that the coefficients of $\chi_\mathcal{A}(q+1)$ form a sign-alternating log-concave sequence of integers with no internal zeros.

\subsection{The variety of critical points}

Finally, the geometry comes into the scene. We are given an essential arrangement $\mathcal{A}$ of $n$ affine hyperplanes in $\mathbb{C}^r$. Our goal is find a subvariety of a product of projective spaces, whose fundamental class encodes the coefficients of the translated characteristic polynomial $\chi_\mathcal{A}(q+1)$. 

The choice of the subvariety is suggested by an observation of Varchenko on the critical points of the master function of an affine hyperplane arrangement \cite{Varchenko}. 
Let $L_1,\ldots,L_n$ be the linear functions defining the hyperplanes of $\mathcal{A}$. A \emph{master function} of $\mathcal{A}$ is a nonvanishing holomorphic function defined on the complement $\mathbb{C}^r \setminus \mathcal{A}$ as the product of powers
\[
\varphi_{\bf u}:=\prod_{i=1}^n L_i^{u_i}, \qquad \mathbf{u}=(u_1,\ldots,u_n) \in \mathbb{Z}^n.
\]

\begin{Var}
If the exponents $u_i$ are sufficiently general, then all critical points of $\varphi_{\bf u}$ are nondegenerate, and
the number of critical points is equal to $(-1)^r \chi_\mathcal{A}(1)$.
\end{Var}

Note that $(-1)^r \chi_\mathcal{A}(1)$ is equal to the number of bounded regions in the complement $\mathbb{R}^r \setminus \mathcal{A}$ when $\mathcal{A}$ is defined over the real numbers, and to the signed topological Euler characteristic of the complement $\mathbb{C}^r \setminus \mathcal{A}$. The conjecture is proved by Varchenko in the real case \cite{Varchenko}, and by Orlik and Terao in general \cite{Orlik-Terao}.

In order to encode all the coefficients of $\chi_\mathcal{A}(q+1)$ in an algebraic variety, we consider the totality of critical points of all possible (multivalued) master functions of $\mathcal{A}$. More precisely, we define the \emph{variety of critical points} $\mathfrak{X}(\mathcal{A})$ as the closure
\[
\mathfrak{X}(\mathcal{A})=\overline{\mathfrak{X}^\circ(\mathcal{A})} \subseteq \mathbb{P}^r \times \mathbb{P}^{n-1}, \qquad \mathfrak{X}^\circ(\mathcal{A})=\Bigg\{\sum_{i=1}^n u_i \cdot \text{dlog}(L_i)(x)=0 \Bigg\} \subseteq (\mathbb{C}^r \setminus \mathcal{A}) \times \mathbb{P}^{n-1},
\]
where $\mathbb{P}^{n-1}$ is the projective space with the homogeneous coordinates $u_1,\ldots,u_n$. 
The variety of critical points first appeared implicitly in \cite{Orlik-Terao2}, and further studied in \cite{Cohen-Denham-Falk-Varchenko,Denham-Garrousian-Schulze}.
See also \cite[Section 2]{HuhML}.

The variety of critical points is irreducible because $\mathfrak{X}^\circ(\mathcal{A})$ is a projective space bundle over the complement $\mathbb{C}^r \setminus \mathcal{A}$.
The cardinality of a general fiber of the second projection 
\[
\text{pr}_2 : \mathfrak{X}(\mathcal{A}) \longrightarrow \mathbb{P}^{n-1}
\] 
is equal to $(-1)^r\chi_\mathcal{A}(1)$, as stated in Varchenko's conjecture. More generally, we have
\[
\big[\mathfrak{X}(\mathcal{A})\big] =
\sum_{i=0}^{r} v_i \big[ \mathbb{P}^{r-i} \times \mathbb{P}^{n-1-r+i}\big] \in H_{2n-2}(\mathbb{P}^r \times \mathbb{P}^{n-1};\mathbb{Z}),
\]
where $v_i$ are the coefficients of the characteristic polynomial
\[
\chi_\mathcal{A}(q+1)=\sum_{i=0}^{r} (-1)^i v_i \hspace{0.7mm} q^{r-i}.
\]
The previous statement is \cite[Corollary 3.11]{HuhML}, which is essentially the geometric formula for the characteristic polynomial of Denham, Garrousian, and Schulze \cite[Theorem 1.1]{Denham-Garrousian-Schulze}, modulo a minor technical difference pointed out in \cite[Remark 2.2]{HuhML}. 
A conceptual proof of the geometric formula can be summarized as follows \cite[Section 3]{HuhML}:

\begin{enumerate}[1.]
\item Applying a logarithmic version of the Poincar\'e-Hopf theorem to a compactification of the complement $\mathbb{C}^r \setminus \mathcal{A}$, one shows that the fundamental class of the variety of critical points captures the characteristic class of $\mathbb{C}^r \setminus \mathcal{A}$.
\item The characteristic class of $\mathbb{C}^r \setminus \mathcal{A}$ agrees with the characteristic polynomial $\chi_\mathcal{A}(q+1)$, because the two are equal at $q=0$ and satisfy the same inclusion-exclusion formula.
\end{enumerate}
See \cite[Section 3]{Denham-Garrousian-Schulze} for a more geometric approach. 

The proof of Theorem \ref{main} is completed by applying Theorem \ref{LC} to the fundamental class of the variety of critical points of $\mathcal{A}$.
\qed

Simple examples show that equalities may hold throughout in the inequalities of Theorem \ref{main}. For example, if $M$ is the uniform matroid of rank $r+1$ with $r+2$ elements, then
\[
h_i\big(\text{IN}(M)\big) = h_i\big(\text{BC}(M)\big)  =1, \qquad i=1,\ldots,r.
\]
However, a glance at the list of $h$-vectors of small matroid complexes generated in \cite{DeLoera-Kemper-Klee} suggests that there are stronger conditions on the $h$-vectors
than those that are known or conjectured. The answer to the interrogative title of \cite{Wilf} seems to be out of reach at the moment.

\end{document}